
\input amstex

\documentstyle{amsppt}

\magnification 1200

\NoBlackBoxes

\define\LL{\llcorner}

\define\join{\times \negthickspace\negthickspace\negthickspace \times}

\define\spt{\operatorname{spt}}

\define\clos{\operatorname{clos}}

\define\graph{\operatorname{graph}}

\define\rank{\operatorname{rank}}

\define\lk{\operatornamewithlimits{lk}}

\loadbold

\define\bull{\boldkey .}

\define\PP{\Bbb P}

\define\CP{\Bbb C \Bbb P}

\define\RP{\Bbb R \Bbb P}

\define\lb{\lbrack\!\lbrack}

\define\rb{\rbrack\!\rbrack}

\topmatter

\title Stiefel-Whitney classes and the conormal cycle of a singular
variety\endtitle

\rightheadtext{Stiefel-Whitney classes}

\author Joseph H. G. Fu and Clint McCrory\endauthor

\leftheadtext{Fu and McCrory}

\date August 21, 1995 \enddate

\address Mathematics Department, University of Georgia, Athens GA 30602
\endaddress

\email fu\@math.uga.edu, clint\@math.uga.edu \endemail

\thanks Research supported in part by NSF grant DMS-9403887. First author
also partially supported by NSF grant DMS-9404366. We thank Adam
Parusi\'nski for his encouragement and help. \endthanks

\abstract A geometric construction of Sullivan's Stiefel-Whitney homology
classes of a real analytic variety $X$ is given by means of the conormal
cycle of an embedding of $X$ in a smooth variety. We prove that the
Stiefel-Whitney classes define additive natural transformations from
certain constructible functions to homology. We also show that, for a
complex analytic variety, these classes are the mod 2 reductions of the
Chern-MacPherson classes.\endabstract



\endtopmatter

\document

We present a new definition of the Stiefel-Whitney homology classes of a
possibly singular real analytic variety $X$. The original definition, due
to Sullivan \cite{S}, involves a triangulation of $X$; its geometric
meaning is unclear. Our definition uses the {\it conormal cycle\/} of an
embedding of $X$ in a smooth variety. The conormal cycle of a subanalytic
subset $X$ of an analytic manifold $M$ was defined by the first author
\cite{F4} using geometric measure theory. The conormal cycle is an
integral current representing (up to sign) Kashiwara's {\it characteristic
cycle\/} of the sheaf $D_M\Bbb R_X$ \cite{F4, 4.7}. Our construction of
the Stiefel-Whitney classes is based on the fundamental observation that
the conormal cycle of a real analytic variety is antipodally symmetric mod
2. We use this observation to give a definition of Stiefel-Whitney classes
which is parallel to the first author's definition of the
Chern-Schwartz-MacPherson classes of a complex analytic variety \cite{F5}.
(This definition of Chern classes is related to earlier work of
Brylinski-Dubson-Kashiwara \cite{BDK} and Sabbah \cite{Sa}.) In fact we
show that the Stiefel-Whitney homology classes of a complex analytic
variety are the mod 2 reductions of the Chern classes. We prove that our
Stiefel-Whitney classes satisfy axioms similar to the Deligne-Grothendieck
axioms for Chern classes, and we prove a specialization formula for the
Stiefel-Whitney classes of a family of varieties. We show that the
Stiefel-Whitney classes of an affine real analytic variety $X$ are
represented by the {\it polar cycles} of $X$, introduced for simplicial
spaces by Banchoff \cite{B} and McCrory \cite{Mc}, and we give a new proof
of the combinatorial formula for Stiefel-Whitney classes of manifolds.

Our central result is a specialization formula for the conormal cycle
(Theorem 3.7). We use this formula to prove the basic results of Kashiwara
and Schapira's calculus of subanalytically constructible functions
\cite{Sc}, \cite{KS, 9.7}, as well as the pushforward and specialization
formulas for Stiefel-Whitney classes.

An advantage of using geometric measure theory to define Chern classes and
Stiefel-Whitney classes is that explicit representative currents can be
constructed. The geometric operations of slicing and limits of subanalytic
currents correspond directly to intersection and specialization of
homology classes. Integral currents have been used by Hardt and McCrory
\cite{HM} to define Steenrod operations, and by Harvey and Zweck \cite{HZ}
to define Stiefel-Whitney classes of vector bundles.

\heading 1. Combinatorial Stiefel-Whitney classes \endheading

In his 1935 thesis, Stiefel defined a characteristic class $w_i(X)$ in the
mod 2 homology of $X$, for every smooth manifold $X$ and each $i\geq 0$:
$w_i(X)$ is the homology class of the singular locus of a general set of
$i+1$ vector fields on $X$ \cite{St}. He conjectured that $w_i(X)$ is
represented by the sum of all the $i$-simplices in the barycentric
subdivision of a triangulation of $X$ \cite{St, p\. 342}. This conjecture
was proved in 1939 by Whitney, who published a sketch of his proof
\cite{W}. Whitney also reinterpreted Stiefel's classes as characteristic
{\it cohomology\/} classes of the tangent bundle of $X$.

Cheeger's rediscovery and proof of the combinatorial formula for
Stiefel-Whitney classes of manifolds in 1969 \cite{Ch} led to Sullivan's
definition of Stiefel-Whitney homology classes of real analytic varieties
with arbitrary singularities. The following result is essentially due to
Sullivan \cite{S}; we update his result using subanalytic technology. (For
a proof that every variety has a subanalytic triangulation, see
\cite{Hi2}, \cite{Ha5}.)

\proclaim{1.1 Theorem} Let $X$ be a compact real analytic variety, let $K$
be a subanalytic triangulation of $X$, and let $K^\prime$ be the
barycentric subdivision of $K$. For each $i\geq 0$, let $s_i(K) \in
C_i(K^\prime;\Bbb Z/2\Bbb Z)$ be the sum of all the $i$-simplices of
$K^\prime$.

$(1)$ The chain $s_i(K)$ is a cycle mod 2.

$(2)$ The homology class $w_i(X) = [s_i(K)] \in H_i(X;\Bbb Z/2\Bbb Z)$ is
independent of the triangulation $K$.

$(3)$ If $f:X\to Y$ is a subanalytic map of real analytic varieties such
that the Euler characteristic $\chi(f^{-1}(y))$ is odd for all $y\in Y$,
then $f_*(w_i(X))=w_i(Y)$.

$(4)$ If $X$ is nonsingular of dimension $d$, then $w_i(X)$ is Poincar\'e
dual to the $(d-i)$th Stiefel-Whitney cohomology class of the tangent
bundle of $X$.

Furthermore, $w_i(X)$ is uniquely characterized by $(3)$ and $(4)$.
\endproclaim

\demo{Proof} Statement (1) of the theorem follows from the fact that $X$
is a (mod 2) Euler space \cite{A}. Several proofs are known that real
analytic varieties are Euler spaces ({\it cf\.} \cite{S}, \cite{BV},
\cite{Ha2}); a new proof is given in (4.4) below. Statements (2) and (3)
follow from Sullivan's mapping cylinder argument \cite{S} ({\it cf\.}
(4.12) below) and the uniqueness of subanalytic triangulations \cite{SY}.
A complete proof of statement (4) (for smooth triangulations) was first
published by Halperin and Toledo \cite{HT}. (The method of Cairns
\cite{C}, applied to a real analytic submanifold of Euclidean space, gives
a triangulation which is both subanalytic and smooth.) The uniqueness of
$w_i(X)$ follows from the representability of mod 2 homology by analytic
manifolds ({\it cf\.} (4.15) below): For every compact real analytic
variety $X$ there exists a subanalytic map $f:V\to X$ such that $V$ is a
compact smooth analytic variety of the same dimension as $X$ and $f_*[V] =
[X]$, where $[X]$ is the mod 2 fundamental class of $X$ \cite{BH}. It
follows that for every compact real analytic variety $X$ there exists a
subanalytic map $\varphi:W\to X$ such that $W$ is a finite disjoint union
of compact smooth analytic varieties (of various dimensions) and
$\chi(\varphi^{-1}(X))$ is odd for all $x\in X$.\qed \enddemo

The chain $s_i(K)$ is the $i$th {\it Stiefel chain\/} of the triangulation
$K$. The homology class $w_i(X)$ is the $i$th {\it Stiefel-Whitney
homology class} of the variety $X$. Clearly $\epsilon\, w_0(X)\equiv
\chi(X)\pmod{2}$, where $\epsilon: H_0(X;\Bbb Z/2\Bbb Z)\to \Bbb Z/2\Bbb
Z$ is the augmentation homomorphism. Also, if $X$ has pure dimension $d$,
then $w_d(X)=[X]$, the mod 2 fundamental class of $X$. However, from the
combinatorial viewpoint the significance of the classes $w_i(X)$ for
$0<i<d$ is not apparent.

\heading 2. The conormal cycle \endheading

Let $X$ be a compact subanalytic subset of the $n$-dimensional oriented
real analytic manifold $M$, and let $f:M\to[0,\infty)$ be a proper,
locally Lipschitz subanalytic function with $f^{-1}(0)=X$. Let
$\pi:T^*M\to M$ be the projection, let $S^*M$ be the cotangent ray space
of $M$, and let $\nu :T^*M-(0)\to S^*M$ be the quotient map.

\definition{2.1 Definition} The {\it conormal cycle} of $X$ in $M$ is the
$(n-1)$-dimensional closed Legendrian integral current $N^*_M(X)$ in
$S^*M$ given by $$N^*_M(X)= \lim_{\epsilon\to 0}\nu_*\langle\lb
df\rb,f\circ\pi,\epsilon\rangle.$$ \enddefinition

This definition is due to Fu \cite{F3} \cite{F4}. Here $\lb df\rb$ denotes
the {\it differential current} of $f$ \cite{F1}. The expression
$\langle\lb df\rb,f\circ\pi,\epsilon\rangle$ denotes the slice of the
current $\lb df\rb$ by the fiber over $\epsilon$ of the function
$f\circ\pi$, and the limit is taken in the flat norm topology. Since  a
real-valued subanalytic function has isolated critical points, the support
of $\langle\lb df\rb,f\circ\pi,\epsilon\rangle$ is contained in $T^*M-(0)$
for $\epsilon$ sufficiently small. Here we use the theory of slicing for
subanalytic integral currents, as developed by Hardt \cite{Ha3}.

The current $N^*_M(X)$ is independent of the choice of function $f$
defining $X$. This follows from a Morse-theoretic characterization of the
local multiplicity of the conormal cycle, which we now describe. Suppose
$X\subset\Bbb R^n$ is compact, and let $N^*(X)$ be the conormal cycle of
$X$ in $\Bbb R^n$. Let $x\in X$, let $\xi:\Bbb R^n\to\Bbb R$ be a nonzero
linear function, and let $[\xi]\in S^*\Bbb R^n$ be the corresponding ray.
For almost all such $\xi$, the multiplicity at $(x,[\xi])$ of the cycle
$N^*(X)$ equals the local change near $x$ of the Euler characteristic of
the super-level sets of $\xi|X$.

More precisely, let $$\iota_X(x,[\xi]) = \lim_{\epsilon\to
0}\lim_{\delta\to 0} \left(\chi\bigl(X\cap B_\epsilon(x)\cap
\xi^{-1}[\xi(x)-h,\infty)\bigr)\bigr|
^{h=+\delta}_{h=-\delta}\right),\tag{2.2}$$ where $B_\epsilon(x)$ is the
closed ball of radius $\epsilon$ centered at $x$. We identify $S^*\Bbb
R^n$ with $\Bbb R^n\times S^{n-1}$, and we let $\kappa$ be the pullback to
$S^*\Bbb R^n$ of the volume form on $S^{n-1}$, normalized to have total
volume $1$. The following theorem is due to Fu \cite{F1, 4.1}.

\proclaim{2.3 Theorem} Let $X$ be a compact subanalytic subset of $\Bbb
R^n$. The conormal cycle $N^*(X)$ is the unique compactly supported
$(n-1)$-dimensional closed Legendrian integral current $T$ on $S^*\Bbb
R^n$ such that, for all smooth functions $\varphi:S^*\Bbb R^n\to \Bbb R$,
$$T(\varphi\cdot\kappa)=\int_{S^{n-1}}\sum_{x\in X}\iota_X(x,\theta)
\varphi(x,\theta)d\theta.\qed $$\endproclaim

Equivalently, let $\rho:S^*\Bbb R^n\to S^{n-1}$ be the projection to the
fiber. Then $$\langle T,\rho,\theta\rangle=\sum_{x\in
X}\iota_X(x,\theta)\lb(x,\theta)\rb.$$ For almost all $\theta\in S^{n-1}$
the function $\iota_X(\cdot,\theta)$ is well-defined and takes nonzero
values at only finitely many points.

\remark{Remarks} (1) It follows from \cite{GM, 3.11} that
$\iota_X(x,[\xi])=1-\chi(\ell^+)$, where $\ell^+$ is the {\it upper
half-link} of $X$ at the point $x$ with respect to the function $\xi$. We
use the super-level sets $\xi^{-1}[a,\infty)$ rather than the sub-level
sets $\xi^{-1}(-\infty,a]$ so that for a convex subset $X$ of $\Bbb R^n$,
the support of the conormal cycle comprises the {\it outward} normals of
$X$ rather than the {\it inward} normals.

(2) For the closed subanalytic subset $X$ of the oriented analytic
manifold $M$, let $CC(\Bbb R_X)$ be Kashiwara's {\it characteristic cycle}
of the sheaf $\Bbb R_X$ on $M$ \cite{KS, 9.4}. Let $a:T^*M\to T^*M$ be the
antipodal map (multiplication by $-1$ in the fiber).
 Now let $\overrightarrow N^*_M(X)$ be the conic Lagrangian cycle of
$T^*M$ corresponding to the conormal cycle $ N^*_M(X)$. Then
$$a_*\overrightarrow N^*_M(X) = CC(\Bbb R_X).$$ For a proof of this
result, see \cite{F4, 4.7}. It follows from \cite{KS, 9.4.4} (or from
(3.12) below) that $$(-1)^n\overrightarrow N^*_M(X)=CC(D_M\Bbb R_X),$$
where $n=\dim M$ and $D_M\Bbb R_X$ is the dual of the sheaf $\Bbb R_X$.
\endremark

If $X$ is a closed subanalytic subset of the analytic manifold $M$, by a
{\it stratification} of $X$ we mean a Whitney stratification with
subanalytic strata \cite{Hi1} \cite{Ha4}. If $\Cal S$ is a stratification
of $X$, then $$\spt N^*_M(X)\subset\bigcup_{S\in\Cal S} n^*S,\tag{2.4}$$
where $n^*S\subset S^*M$ is the conormal ray bundle of the stratum $S$ in
$M$ \cite{F4, 4.6.6}.

It will be convenient to extend the definition of the conormal cycle to
all subanalytic sets. If $X$ is a subanalytic subset of the oriented real
analytic manifold $M$, then $X$ is locally closed in $M$, and the current
$N^*_M(X)$ is uniquely determined by the following property. If $U$ is a
subanalytic open subset of $M$ such that the closure $\overline U$ is
compact and $X\cap \overline U$ is closed in $\overline U$, then
$$N^*_M(X)\ \LL\ \pi^{-1}U = N^*_M(X\cap \overline U)\ \LL\
\pi^{-1}U,\tag{2.5}$$ and $N^*_M(X)\ \LL\ \pi^{-1}U=0$ if $X\cap \overline
U=\emptyset$.

The following theorem of Fu \cite{F4, 4.6} is of fundamental importance
for applications of the conormal cycle.

\proclaim{2.6 Theorem} Let $X$ be a closed stratified subanalytic subset
of the oriented analytic $n$-manifold $M$, with $\dim X < n$. Let
$f:M\to\Bbb R$ be an analytic function such that for every stratum $S$ of
$X$, $0$ is a regular value of $f|S$. Then $$N^*_M\left(X\cap
f^{-1}(-\infty,0]\right)=N^*_M\left(X\cap f^{-1}(-\infty,0)\right) +
j_*\left(\partial N^*_M\left(X\cap
f^{-1}(-\infty,0)\right)\times[0,1]\right),$$ where $j:S^*M\times\Bbb R\to
S^*M$ is the join map $$j(\xi,t)=[(1-t)\xi+tdf_{\pi(\xi)}].$$ In fact,
$N^*_M(X\cap f^{-1}(-\infty,0])$ is the unique Legendrian cycle $T$ such
that \roster \item"(i)" $\spt T \subset \pi^{-1}(X\cap f^{-1}(-\infty,
0])$, \item"(ii)" $T\ \LL\ \pi^{-1}(f^{-1}(-\infty,0))= N^*_M(X\cap
f^{-1}(-\infty,0))$, and \item"(iii)" if $f(x)=0$ then $[-df_x]\notin\spt
T$.\qed \endroster\endproclaim

Note that by replacing $f$ by $-f$ we have $$N^*_M\left(X\cap
f^{-1}[0,\infty)\right)=N^*_M\left(X\cap f^{-1}(0,\infty)\right)+ k_*\left
(\partial N^*_M\left(X\cap f^{-1}(0,\infty)\right)\times[0,1]\right),$$
$$k(\xi,t)=[(1-t)\xi-tdf_{\pi(\xi)}].$$

We make extensive use of the following {\it product formula} for conormal
cycles. Let $\overrightarrow N^*_M(X)$ be the conic Lagrangian cycle of
$T^*M$ corresponding to the conormal cycle $ N^*_M(X)$. Then for closed
subanalytic sets $X\subset M$, $X'\subset M'$, we have $$\overrightarrow
N^*_{M\times M'}(X\times X') = \overrightarrow N^*_M(X)\times
\overrightarrow N^*_{M'}(X')\tag{2.7}$$ \cite{F4, 4.5}. In terms of the
conormal cycle in the ray space, this becomes $$\align N^*_{M\times
M'}(X\times X') &= N^*_M(X) \join N^*_{M'}(X')\tag{2.8}\\ &=
j_*\left(N^*_M(X)\times [0,1]\times N^*_{M'}(X')\right), \endalign$$ where
$j([\xi],t,[\eta])= [t\xi+(1-t)\eta]\in S^*(M\times M')$.

The relation of the conormal cycle to the Euler characteristic is
reflected in two basic properties. First, the conormal cycle is {\it
additive} \cite{F4, 4.2}: If $X$ and $Y$ are closed subanalytic subsets of
the oriented analytic manifold $M$, then $$ N_M^*(X\cup Y)=
 N_M^*(X)+  N_M^*(Y)-  N_M^*(X\cap Y).\tag{2.9}$$ Second, the conormal
cycle satisfies a {\it Gauss-Bonnet property} \cite{F4, 1.5}: For $X$ a
compact subset of $M$ as above, if $M$ is Riemannian, with
Chern-Gauss-Bonnet form $\Omega\in\bigwedge^n(M)$ and transgression form
$\Pi\in \bigwedge^{n-1}(S^*M)$, then $$ N_M^*(X)(\Pi)+\int_X\Omega =
\chi(X).\tag{2.10}$$
 If $M=\Bbb R^n$, then the form $\kappa$ on $S^*M=\Bbb R^n\times S^{n-1}$
(2.3) is a (closed) transgression form. Thus if we choose $\theta\in
S^{n-1}$ so that $\Bbb R^n\times\theta$ is transverse to a stratification
of $\spt N^*_{\Bbb R^n}(X)$, then $$( N^*_{\Bbb R^n}(X) \cdot\lb\Bbb
R^n\times\theta\rb)(1)=\chi(X).\tag{2.11}$$ Here we use the notation
$A\cdot B$ to denote the intersection of the integral currents $A$ and
$B$, and the notation $\lb Z\rb$ to denote the $k$-dimensional integral
current associated to the oriented $k$-manifold $Z$.

As a consequence of (2.6) and the Gauss-Bonnet property (2.10) we have the
following result on the Morse theory of the conormal cycle.

\proclaim{2.12 Proposition} Let $X$ be a compact subanalytic subset of the
oriented analytic $n$-manifold $M$, with $\dim X < n$, and let $f:M\to\Bbb
R$ be an analytic function. If $df_x\neq 0$ and $[-df_x]\notin\spt
N^*_M(X)$ for all $x\in f^{-1}[a,b]$, and if $$N^*(X)\ \LL\
\pi^{-1}\left(f^{-1}(a)\right) = 0, \tag{1} $$then $$\chi\left(X\cap
f^{-1}(-\infty,a]\right)= \chi\left(X\cap
f^{-1}(-\infty,b]\right).$$\endproclaim

\demo{Proof} The hypothesis implies that the join map $j$ of (2.6) is
well-defined and smooth on $\spt N^*(X) \cap
\pi^{-1}\left(f^{-1}[a,b]\right)$. Let $T=j_*\left( N^*\left(X\cap
f^{-1}(a,b)\right) \times [0,1]\right)$, an $n$-dimensional integral
current. By (2.6), $$\align  &N^*\left(X \cap f^{-1}(-\infty ,b]\right) -
  N^*\left(X \cap f^{-1}(-\infty,a]\right)\\ &=  N^*\left(X \cap
f^{-1}(a,b)\right) + j_*\left(\partial \left( N^*\left(X \cap
f^{-1}(-\infty,b)\right) -   N^* \left(X \cap
f^{-1}(-\infty,a)\right)\right)\times [0,1]\right) \\ &=\partial
T,\endalign$$ since, by (1), $\partial\left( N^*(X)\ \LL\
\pi^{-1}\left(f^{-1}(-\infty ,a)\right)\right) = -\partial\left( N^*(X)\
\LL\ \pi^{-1}\left(f^{-1}(a,\infty )\right)\right)$. Therefore,
 $$ \align \chi\left( X \cap f^{-1}(-\infty,b]\right) &- \chi\left( X \cap
f^{-1}(-\infty,a]\right)\\ &= \left( N^*\left(X \cap
f^{-1}(-\infty,b]\right) -
 N^* \left(X \cap f^{-1}(-\infty,a]\right)\right)(\Pi)\
 \ \text { by (2.10) } \\   &= \partial T (\Pi) \\ &= T (d \Pi) \\ &= T
(\pi^*(\Omega)) \\ &= \pi_*T(\Omega) \\ &= 0, \endalign$$ since $\dim \spt
\pi_* T < n$. \qed \enddemo

\subhead 2.13 The mod 2 conormal cycle\endsubhead Let $X$ be a compact
subanalytic subset of the real analytic $n$-manifold $M$. If we do not
assume that $M$ is oriented, the conormal cycle of $X$ in $M$ can be
defined as an integral current modulo 2. It is the $(n-1)$-dimensional
closed Legendrian integral current mod 2 given by $$\Cal
N^*_M(X)=\lim_{t\to 0}\nu_*\langle\lb f\rb,f\circ\pi,\epsilon\rangle$$ as
in (2.1), where $\lb f\rb$ is the mod 2 differential current of a
non-negative function $f$ which defines $X$, and the slice and limit are
in the sense of subanalytic integral currents mod 2 ({\it cf\.}
\cite{Ha1}, \cite{Ha3}).

Almost all the results of this section are true for the mod 2 conormal
cycle; only the Gauss-Bonnet property (2.10) requires that the ambient
manifold $M$ is oriented. The proofs of the mod 2 results are parallel to
the proofs in the oriented case. The main technical point is the mod 2
version of Fu's uniqueness theorem for Lagrangian cycles \cite{F1, 1.1},
which is used to characterize the mod 2 differential current $\lb f\rb$
and to prove the uniqueness part of (2.6).

An alternate approach is to define the mod 2 conormal cycle in terms of
local orientations of $M$. If $\{U_i\}$ is a covering of $M$ by oriented
open sets, then $\Cal N^*_M(X)$ is uniquely characterized by the property
that $$\Cal N^*_M(X)\ \LL\ \pi^{-1}U_i = \{N^*_{U_i}(X\cap U_i)\}_2$$ for
all $i$, where $\{C\}_2$ denotes the mod 2 reduction of the integral
current $C$. Thus the local properties of $\Cal N^*$ can be proved by
reduction mod 2 of the corresponding properties of $N^*$.

\heading 3. Specialization\endheading

As Parusi\' nski has emphasized, the process of {\it specialization} is of
basic importance in subanalytic homology theory. In this section we prove
a specialization formula for the conormal cycle, and we apply it to derive
those aspects of Kashiwara and Schapira's calculus of subanalytically
constructible functions \cite{Sc} \cite{KS, 9.7} which are needed to state
and prove the Deligne-Grothendieck axioms for Stiefel-Whitney classes. (In
contrast, Kashiwara and Schapira define operations on constructible
functions using standard operations of the derived category of sheaves.)
We define the operations of Euler characteristic (integral), pushforward
(direct image) and duality, and we prove the key relations between these
operations.

\definition{3.1 Definition} Let $M$ be a real analytic manifold. The
function $\alpha:M\to \Bbb Z$ is (subanalytically) {\it constructible} if
$\alpha^{-1}(n)$ is a subanalytic set for all $n\in\Bbb Z$, and the
collection $\{\alpha^{-1}(n)\}_{n\in\Bbb Z}$ is locally finite.
\enddefinition

The set of constructible functions on $M$ is a ring, with
$(\alpha+\beta)(p)=\alpha(p)+\beta(p)$ and $(\alpha\cdot\beta)(p) =
\alpha(p)\beta(p)$. Every constructible function on $M$ can be written as
a locally finite sum $$\alpha= \sum_in_i1_{X_i},\tag{3.2}$$ where $n_i\in
\Bbb Z$, each $X_i$ is a closed subanalytic subset of $M$, and $1_X$
denotes the characteristic function of $X$.

\definition{3.3 Definition} If $\alpha$ is a constructible function on $M$
with compact support, the {\it Euler characteristic} $\chi(\alpha)\in\Bbb
Z$ is defined as follows. Write $\alpha$ as a finite sum (3.2) with each
$X_i$ compact. Then $$\chi(\alpha)=\sum_in_i\chi(X_i).$$\enddefinition

Kashiwara and Schapira use the suggestive notation
$\chi(\alpha)=\int_M\alpha.$ By the triangulation theorem for subanalytic
sets \cite{Hi2} \cite{Ha5}, the Euler characteristic is uniquely
determined by two properties: \roster\item $\chi(\alpha+\beta)=
\chi(\alpha)+\chi(\beta)$ for all $\alpha$, $\beta$. \item $\chi(1_X)=1$
if $X$ is subanalytically homeomorphic to a closed ball.\endroster

\definition{3.4 Definition} If $\alpha$ is a constructible function on the
oriented analytic manifold $M$, the {\it conormal cycle} $N^*_M(\alpha)$
is defined as follows. Write $\alpha$ as a locally finite sum (3.2), and
let $$N^*_M(\alpha) = \sum_in_iN^*_M(X_i),$$ a closed Legendrian integral
current in the cotangent ray space $S^*M$. \enddefinition

By (2.9) the conormal cycle $N^*_M(\alpha)$ is well-defined. The
Gauss-Bonnet theorem (2.10) yields $$N_M^*(\alpha) (\Pi) = \chi(\alpha)
\tag{3.5}$$ for a compactly supported constructible function $\alpha$ with
$\dim \spt \alpha < \dim M$, where $\Pi$ is a transgression form for $M$.

\proclaim{3.6 Proposition} If $M$ is an oriented analytic manifold, then
$N^*_M$ is an isomorphism from the group of constructible functions
$\alpha$ on $M$ such that $\spt\alpha$ has codimension at least 2 to the
group of subanalytic Legendrian cycles $Z$ such that $\pi(\spt Z)$ has
codimension at least 2. \endproclaim

\demo{Proof} Suppose $\alpha,\beta$ are constructible functions with
support of codimension at least 2, with $ N^*_M(\alpha)= N^*_M(\beta)$ and
$\dim(\spt\alpha)\geq \dim(\spt\beta)$. Let $\Cal S$ be a subanalytic
Whitney stratification of $M$ so that both $\alpha$ and $\beta$ are
constant on the strata of $\Cal S$. Let $S$ be a $k$-dimensional stratum
of $\Cal S$, where $k=\dim(\spt\alpha)$. Then $ N^*_M(\alpha)\ \LL\
\pi^{-1}(S) = m N^*_M(S)$ and $ N^*_M(\beta)\ \LL\ \pi^{-1}(S) = n
N^*_M(S)$ for some $m$, $n\in\Bbb Z$. Thus $m=n$, and we obtain
$\alpha=\beta$ by induction on $\dim(\spt\alpha)$.

On the other hand, suppose $Z$ is subanalytic Legendrian cycle  in $S^*M$,
put $A= \spt Z$, and suppose that the codimension $k$ of $\pi(A)$ is at
least 2. Let $\Cal S$ be a stratification of $A$ such that the restriction
of $\pi $ to each $S \in \Cal S$ is a submersion. The Legendre condition
implies that each $S \subset \nu^*\pi(S)$. Let
 $\Cal V = \{ \pi(S) - \bigcup_{S' \in \Cal S, S' \ne S} \clos \pi(S')
 \mid \pi(S) \text{ has maximal dimension }m-k \} $, where $m = \dim M$;
this is a collection of nonempty open subanalytic submanifolds $V$ of $M$,
and each $V \in \Cal V $ has a neighborhood $U \subset M$ such that $\spt
(Z \ \LL \ \pi^{-1}(U)) \subset \nu^* V$. Since $ k \ge 2$, each $\nu^* V$
is connected, so the constancy theorem implies that
 $Z \ \LL \ \pi^{-1}(U) = n_V \lb \nu^*V\rb$ for some $n_V \in \Bbb Z$.
Then $\pi\left(\spt \left(Z - \sum n_V N^*_M\left(\overline
V\right)\right)\right)$ has codimension larger than $k$, and the assertion
follows by induction. \qed\enddemo

The central result of this section is the following specialization formula
for the conormal cycle. The following notation will be useful. If
$g:M\to\Bbb R$, $X\subset M$, and $s$, $t\in\Bbb R$, $s<t$, let $X_t=X\cap
g^{-1}(t)$ and $X_{[s,t]}=X\cap g^{-1}([s,t]).$ For $p$, $q\in M$, let
$\rho_p(q)$ denote the distance from $p$ to $q$ in an analytic metric on
$M$.

\proclaim{3.7 Theorem} Let $\alpha$ be a constructible function on the
oriented analytic manifold $M$ such that $X=\spt\alpha$ is compact and has
positive codimension. Let $g:M \to \Bbb R$ be a subanalytic function. For
$t\in \Bbb R$, put $\alpha_t = \alpha\cdot 1_{X_t}$. For each $p \in X_0$,
there is $\epsilon_0 > 0$ such that if $0 < \epsilon < \epsilon_0$ then
the limit $$\phi_0\alpha(p) = \lim_{t\downarrow 0} \chi(\alpha_t\cdot
1_{B_\epsilon(p)})$$ exists, and this limit is independent of $\epsilon$.
The resulting function $\phi_0\alpha$ is constructible, and $$
N^*_M(\phi_0\alpha)=\lim_{t\downarrow 0} N^*_M(\alpha_t),$$ where the
limit is in the flat norm topology. More precisely, for all sufficiently
small $h>0$ there is an integral current $Z_h$ in $S^*M$ such that \roster
\item "{(a)}"$\partial Z_h= N^*_M(\alpha_h) - N^*_M(\phi_0\alpha)$, \item
"{(b)}"$\spt Z_h \subset \pi^{-1}(X_{[0,h]})$, \item "{(c)}"$Z_h \to 0$ in
mass as $h \to 0$.\endroster \endproclaim

\demo{Proof} For simplicity we let $N^*(\alpha)=N^*_M(\alpha)$. Let $\Cal
S$ be a Whitney stratification of $X$ such that $\alpha|S$ is constant and
$g|S$ is smooth for all $S \in \Cal S$. Let $h_0 > 0$ be so small that for
all $S\in\Cal S$ the restriction $g|S$ has no critical values in
$(0,h_0]$. Put $$Z= \left(j_* + k_*\right) \left(\left(N^*(\alpha)\
\llcorner\ (g\circ \pi)^{-1}(0,h_0]\right) \times [0,1]\right),$$ where
$j$, $k:S^*M\times\Bbb R\to S^*M$ are the join maps defined in (2.6). Then
$Z$ is a subanalytic integral current of dimension $n=\dim M$, with
$$\align \langle Z, g \circ \pi, t \rangle &= (j_* + k_*)\left(\langle
N^*(\alpha) , g \circ \pi, t \rangle \times [0,1]\right) \\ &=
N^*(\alpha_t) \endalign $$ for $t \in (0,h_0]$, by (2.6) and the addition
formula (2.9). In particular, $$\lim_{t \downarrow 0} N^*(\alpha_t) =
N^*(\alpha_{h_0})-\partial Z$$ exists, and by (3.6) equals $N^*(\psi)$ for
some constructible function $\psi$. To complete the proof we must show
that $ \psi =\phi_0\alpha $; then the currents $Z_h = Z \ \llcorner\
\pi^{-1}g^{-1} (0,h]$, $0 < h \le h_0$, satisfy the stated relations.

 Let $\Cal T$ be a stratification of $X$ subordinate to $\Cal S \cup
\{X_0\}$, such that $\Cal T$ satisfies the $(a_g)$ condition ({\it cf\.}
\cite{Be, Prop\. 10}). Then $$\align
 \spt N^*(\psi) &\subset \limsup_{h \downarrow 0} \left(\spt
N^*(\alpha_{h})\right) \\ &\subset \limsup_{h \downarrow 0}  \bigcup_{S
\in \Cal S} n^*(S_{h}) \\ &\subset \limsup_{h \downarrow 0}  \bigcup_{T
\in \Cal T} n^*(T_{h}) \tag{1}\\ &\subset \bigcup_{T\in \Cal T, \, T
\subset X_0} n^* T , \endalign $$ where the last inclusion follows from
the $(a_g)$ condition. Now write $\psi = \sum n_i 1_{Y_i}$, $n_i\neq 0$,
where the sum is locally finite, and $Y_i$ are compact and subanalytic,
with $Y_i=\overline{Y_i^\circ}$, $Y_i^\circ$ smooth and $\bigcup_i n^*
(Y_i^\circ) \supset \spt N^*(\psi)$, let $p \in X_0$ be fixed, and choose
$ \epsilon_0 > 0$ so small that the following two conditions hold: \roster
\item"{(i)}" For all $r \in (0, \epsilon_0)$, $Y_i \cap B_r(p)$ is
homeomorphic to a cone. (Such $\epsilon_0$ exists by \cite{Ha6}.)
\item"{(ii)}" There exists $c>0$ with $\Vert d_x(\rho_p | T) \Vert \ge c$
for all strata $T \in \Cal T, \, T \subset X_0$, and all points $x \in T
\cap B_{\epsilon_0}(p)$, $x\neq p$. \endroster Then, for $r<\epsilon_0$,
$$\psi(p) = \sum n_i 1_{Y_i} (p) =\sum n_i \chi(Y_i \cap B_r(p)),
\tag{2}$$ and (1) implies that $$\Vert d_x(\rho_p | T_{h}) \Vert \ge
c/2\tag{3}$$ whenever $h >0$ is sufficiently small, $T \in \Cal T$, and $x
\in T_{h}\cap B_r(p)$.

We claim that, if $r<\epsilon_0$, then $$\lim_{h \downarrow 0}
N^*(\alpha_{h}\cdot 1_{B_r(p)}) = N^*(\psi\cdot 1_ {B_r(p)}).\tag{4}$$ For
the restrictions of the two sides to the interior of $B_r(p)$ are equal by
construction. On the other hand, (1), (3) and the join construction of
(2.6) imply that $$\clos\left(\bigcup_{h \le h_0} \spt
N^*\left(\alpha_h\cdot 1_{B_r(p)}\right)\right) \cap \graph (-d\rho_p)
\setminus \pi^{-1}(p) = \emptyset,$$ and therefore $$\spt
\left(\lim_{h\downarrow 0}   N^*\left(\alpha_{h}\cdot
1_{B_r(p)}\right)\right)\cap \graph (-d\rho_p)
 \cap \pi^{-1} S_r(p) = \emptyset.$$ Now (4) follows from the uniqueness
part of (2.6).

Therefore, if $\Pi$ is a transgression form for $M$, we have by (2) and
(2.10) that $$ \align \phi_0\alpha(p)&=\lim_{h \downarrow 0}
\chi(\alpha_{h}\cdot 1_{B_r(p)})\\ &=\lim_{h\downarrow 0}
N^*(\alpha_{h}\cdot 1_{B_r(p)}) (\Pi) \\ &=N^*(\psi\cdot 1_ {B_r(p)})
(\Pi) \\ &= \psi(p), \endalign $$ as desired.\qed\enddemo

\definition{3.8 Definition} Let $M$ and $N$ be analytic manifolds, and let
$f:X\to Y$ be a subanalytic map, where $X\subset M$ and $Y\subset N$. Let
$\alpha$ be a constructible function on $M$ such that $\spt\alpha\subset
X$, and suppose that $f|\spt\alpha$ is proper. The {\it pushforward}
$f_*\alpha$ is the constructible function on $N$, with $\spt
f_*\alpha\subset Y$,  given by $$(f_*\alpha)(q)=\chi\left(\alpha\cdot
1_{f^{-1}(q)}\right).$$ \enddefinition

The properties of pushforward  may be derived from the specialization
formula, applied to the mapping cylinder. Given subanalytic subsets
$X\subset \Bbb R^m, Y \subset
 \Bbb R^n$ and a subanalytic map $f: X \to Y$, consider the {\it mapping
cylinder} of $f$, $$C = \left\{ \left((1-t)x, t, f(x)\right)\  |\ x \in
X,\ 0 \le t \le 1\right\}\ \cup\ \left\{(0,1,y)\ |\ y\in Y\right\}.$$ This
is a subanalytic subset of $\Bbb R^m \times \Bbb R \times \Bbb R^n$ of
positive codimension. Let $g: \Bbb R^m \times \Bbb R \times \Bbb R^n \to
\Bbb R$ be the projection to the middle factor, and let $C_t=C\cap
g^{-1}(t)$. Suppose now that $\alpha$ is a constructible function on $\Bbb
R^m$ with compact support contained in $X$. We let $\tilde\alpha$ be the
constructible function on $\Bbb R^m \times \Bbb R \times \Bbb R^n$ given
by $\tilde\alpha((1-t)x,t,f(x))=\alpha(x)$ for $x \in X$, $t\in [0,1)$,
and $\tilde\alpha$ equal to zero elsewhere. For $t\in[0,1)$, let
$\tilde\alpha_t=\tilde\alpha\cdot 1_{C_t}$. For $s\in[0,1]$, we define the
right and left specializations of $\tilde\alpha_t$ as $t\to s$ to be the
constructible functions $\phi_{s}^\pm \tilde\alpha$ given by $\phi_{s}^\pm
\tilde\alpha(p)=\lim_{t\to s\pm} \chi(\tilde \alpha_t\cdot
1_{B_\epsilon(p)})$ for sufficiently small $\epsilon>0$.

\proclaim{3.9 Proposition} Let $X\subset\Bbb R^m$, $Y\subset\Bbb R^n$ be
subanalytic sets. Let $f:X\to Y$ be a subanalytic map with mapping
cylinder $C$, and let $\alpha$ be a constructible function on $\Bbb R^m$
with compact support contained in $X$. Then $$\align\phi_0^+\tilde\alpha
&= \alpha, \tag{1}\\ \phi_1^-\tilde\alpha &= f_*\alpha. \tag{2}
\endalign$$ Furthermore, there is a subanalytic integral current $W$ with
$\spt W \subset \pi^{-1}C$ and $$\partial W = N^*_M(f_*\alpha) -
N^*_M(\alpha).\tag{3}$$ Here we identify the constructible functions
$\alpha$ and $f_* \alpha$ with the corresponding functions on $\Bbb R^m
\times \{0\} \times \{0\}$ and $\{0\} \times\{0\} \times \Bbb R^n$
respectively. \endproclaim

\demo{Proof}  First we establish the specialization formulas.

We prove (2) first. Let $\Cal S$ be a Whitney stratification of $C$ and
$\Cal T$ an $(a_g)$ stratification refining $\Cal S \cup C_1$, as in the
proof of (3.7); let $q \in C_1=Y$ be a given point. Let $\epsilon_0 > 0$
be so small that $\Vert  d_x (\rho_q | T) \Vert \ge c_0 > 0$ for $x \in T
\cap B_{\epsilon_0}(q), \, x \ne q$, $C_1 \supset T \in \Cal T$. Theorem
2.6 and the analogue of the relation (3.7)(1) imply that $$[-d_x\rho_q]
\notin \limsup_{t \uparrow 1} \left(\spt N^*\left(\tilde\alpha_t \cdot 1_
{B_{\epsilon_0}(q)}\right)\right) \tag{4}$$ for $x \in S_{\epsilon_0}(q)$.
 Consider the modified distance function $\tilde \rho (x,t,y) =
\rho_q(y)$; clearly $\rho_q = \tilde\rho$ on $C_1$. Thus (4) and the
$(a_g)$ condition imply that, given $\epsilon_1 \in(0,\epsilon_0)$, we
have for all $t$ sufficiently close to $1$,
 $$[-d_x \tilde \rho] \notin \spt N^*\left(\tilde\alpha_t \cdot 1_
{B_{\epsilon_0}(q)}\right)$$ for $x \in
\tilde\rho^{-1}[\epsilon_1,\epsilon_0]$, and,
 if $\epsilon_1$ is sufficiently small, $C_t\cap\tilde\rho^{-1}[0,
\epsilon_0]\supset C_t\cap B_{\epsilon_0}(q) \supset
C_t\cap\tilde\rho^{-1}[0,\epsilon_1].$ Therefore for such $t$ (2.12)
yields $$\align \chi(\tilde\alpha_t\cdot 1_{B_{\epsilon_0}(q)})&=
\chi\left(\tilde\alpha_t \cdot 1_{\tilde\rho^{-1}[0,\epsilon_1]}\right)\\
&=\chi\left(\alpha\cdot 1_{f^{-1}\left(B_{\epsilon_1}(q)\right)}\right)\\
&=\chi\left(\alpha\cdot 1_{f^{-1}(q)}\right),\endalign$$ since $C_t \cap
\tilde\rho^{-1}[0,\epsilon_1]\cong f^{-1}(B_{\epsilon_1}(q)).$ The
relation (2) now follows from the specialization formula (3.7).

To prove (1) we simply repeat the argument of the last paragraph with
$\tilde \rho$ replaced by the function $\hat \rho(x,t,y) = \rho_p(x)$,
where $p \in X$, and we observe that $C_t \cap \hat
\rho^{-1}[0,\epsilon_1] \cong X \cap B_{\epsilon_1}(p)$ for $t$ close to
0.

Now by \cite{Hi1, Lemma 4.8.3} and \cite{Be, Lemma 15}, given a
subanalytic Whitney stratification $\Cal X$ of $X$, there exists a
subanalytic Whitney stratification $\Cal T$ of the mapping cylinder $C$
satisfying $(a_g)$ and with the property that $\Cal T|C_{[0,1)}=\Cal
X\times [0,1).$ More precisely, $C_{[0,1)} = \{(x,t,y)\in C\ |\ 0\leq
t<1\}$ is a union of strata of $\Cal T$, and all the strata of $C_{[0,1)}$
are of the form $H(S\times\{0\})$ or $H(S\times(0,1))$ for $S\in\Cal X$
and $H:X\times[0,1)\to C_{[0,1)}$ the homeomorphism $H(x,t)=(x,t,f(x))$.
Thus for all $T\in\Cal T$, the restriction $g|T$ has no critical values in
$(0,1)$.

Therefore, by the proof of (3.7), there exist integral currents $Z^+$,
$Z^-$ in $S^*(\Bbb R^m\times\Bbb R\times\Bbb R^n)$ such that
$$\align\partial Z^+_0&=N^*(\tilde\alpha_{\frac{1}{2}})-
N^*(\phi^+_0\tilde\alpha),\ \ \spt Z^+_0\subset\pi^{-1}(C_{[0,\frac
{1}{2}]}),\\ \partial Z^-_1&=N^*(\tilde\alpha_{\frac{1}{2}})-
N^*(\phi^-_1\tilde\alpha),\ \ \spt Z^-_1\subset\pi^{-1}(C_{[\frac
{1}{2},1]}).\endalign$$ Thus for $W=Z^+_1-Z^-_0$, (3) follows from (1) and
(2). \qed \enddemo

\proclaim{3.10 Proposition} Suppose $f :X \to Y$ and $g: Y \to Z$ are
proper subanalytic maps. Then the pushforward operations under these maps
satisfy the relation $$(g \circ f)_* = g_* \circ f_*. \tag{1}$$ If the
constructible function $\alpha $ has compact support in $X$ then
$$\chi(f_* \alpha) = \chi(\alpha).\tag{2}$$\endproclaim

\demo{Proof} To prove (2), we imbed $X,Y$ in $\Bbb R^m, \Bbb R^n$
respectively, and then apply (3.9)(3) and the Gauss-Bonnet theorem (2.10).
Using (2), we deduce (1) as follows. It suffices to prove that
$$g_*(f_*\alpha)=(g\circ f)_* \alpha$$ for $X$ compact and $\alpha=1_X$.
Let $\beta=f_*\alpha$ and $\gamma=(g\circ f)_*\alpha$. We will show that
$g_*\beta=\gamma$. Let $z\in Z$, let $X'=(g\circ f)^{-1}(z)$,
$Y=g^{-1}(z)$, and let $f':X'\to Y'$ be the restriction of $f$. Then,
since $f'_*(1_{X'})=\beta\cdot 1_{Y'}$, $$\align \gamma(z) &=
\chi(1_{X'})\\ &= \chi(f'_*(1_{X'}))\ \ \ \ \text{by (2)}\\ &=
\chi(\beta\cdot 1_{Y'})\\ &= (g_*\beta)(z),\endalign$$ as desired.\qed
\enddemo

\remark{Remark} If $f:M\to N$ is a subanalytic map and $\beta$ is a
constructible function on $N$, the {\it pullback} $f^*\beta$ is the
constructible function on $M$ defined by $$(f^*\beta)(p)=\beta(f(p)).$$ It
is easy to see that $(f\circ g)^*=g^*\circ f^*$, but $\chi(f^*\beta)$ is
not equal to $\chi(\beta)$, in general.\endremark

\definition{3.11 Definition} Let $\alpha$ be a constructible function on
the analytic manifold $M$. The {\it dual} of $\alpha$ is the constructible
function $D_M\alpha$ on $M$ defined by $$(D_M\alpha)(p) =\chi\left(\alpha
\cdot 1_{B_\epsilon^\circ(p)}\right) = \alpha(p) - \chi(\alpha \cdot
1_{S_\epsilon(p)}),$$ where $B_\epsilon^\circ(p)$ is the open ball of
radius $\epsilon>0$ about $p$ in an analytic metric on $M$, and $\epsilon$
is sufficiently small. \enddefinition

To see that $D_M\alpha$ does not depend on $\epsilon$ or on the choice of
metric, note that if $\alpha = \sum_i n_i 1_{X_i}$ as in (3.2), then
$D_M\alpha = \sum_i n_i D_M1_{X_i}$. If $X$ is a closed subanalytic subset
of $M$ and $p\in X$, then $$(D_M1_X)(p) = \chi_p(X),$$ where $\chi_p(X)$
is the local Euler characteristic of $X$ at $p$, $\chi_p(X) = \sum_i\rank
H_i(X, X\setminus \{p\}).$ Furthermore, $\chi_p(X) = 1 - \chi(L)$, where
$L$ is a subanalytic link of $p$ in $X$ ({\it cf\.} \cite{CK}).

\proclaim{3.12 Theorem} If $\alpha$ is a constructible function on the
oriented analytic $n$-manifold $M$, then $$(-1)^na_*N^*_M(\alpha) =
N^*(D_M \alpha).$$ \endproclaim

\demo{Proof} Embedding $M$ in a Euclidean space, we may assume by the
product formula (2.7) that $M = \Bbb R^n$ and $\dim\spt\alpha \leq n-2$.
Let $N^*=N^*_{\Bbb R^n}$ and $D=D_{\Bbb R^n}$. Using the isomorphism (3.6)
between Legendrian cycles on $S^*\Bbb R^n$ and constructible functions on
$\Bbb R^n$, there is a linear operator $A$ such that $(-1)^na_*N^*(\alpha)
= N^*(A \alpha)$; we must show that $A = D$. By the triangulation theorem
for subanalytic sets it is enough to prove that $A(1_\Delta) =
D(1_\Delta)$ for every subanalytic simplex $\Delta \subset \Bbb R^n$.

Let $k = \dim \Delta$, and let $f: B \to \Delta$ be a subanalytic
homeomorphism from a $k$-ball $B$ to $\Delta$. Then $$f_*(1_B)=1_\Delta,\
\ \ \ f_*(1_{\partial B}) = 1_{\partial\Delta},$$
$$D(1_B)=(-1)^k(1_B-1_{\partial B}),\ \ \ \
D(1_\Delta)=(-1)^k(1_\Delta-1_{\partial \Delta}),$$ and therefore
$$f_*D(1_B)=D(1_\Delta).$$ Taking $n$ large enough we may assume that $B$
and $\Delta$ are embedded in $\Bbb R^n = \Bbb R^m \times \Bbb R \times
\Bbb R^{n - m-1}$ as the ends of the mapping cylinder $C$ of $f$. Consider
the constructible function $\widetilde{1_{B}}=1_C-(1_{C_0}+1_{C_1})$ on
$C$. Then by (3.9) $$\align \phi_0^+ (\widetilde {1_B})& = 1_B, \tag{1}\\
\phi_1^- (\widetilde {1_B}) &= f_* (1_B)=1_\Delta.\tag{2} \endalign $$

Now we claim that $$A(\widetilde{1_B})_t = D(\widetilde{1_B})_t, \ \ 0 < t
< 1.\tag{3}$$ If $0 < s < t < 1$, then the map $\psi_{s,t}:\Bbb R^n\to\Bbb
R^n$, $$\psi_{s,t}(x,u,y) = (\textstyle\frac{s}{t} x, s+t -u
,\textstyle\frac{1-s}{1-t} y),$$ is an analytic isomorphism taking
$(\widetilde{1_B})_t$ to $(\widetilde{1_B})_s$. Let $\Psi_{s,t}:S^*\Bbb
R^n\to S^*\Bbb R^n$ be the isomorphism induced by $\psi_{s,t}$. Then
$$\align N^*A(\widetilde{1_B})_s &= (-1)^na_*N^*(\widetilde{1_B})_s \\ &=
(-1)^na_*(\Psi_{s,t}^{-1})_*N^*(\widetilde{1_B})_t \\ &=
(\Psi_{s,t}^{-1})_* (-1)^na_* N^*(\widetilde{1_B})_t \\ & =
(\Psi_{s,t}^{-1})_* N^*A (\widetilde{1_B})_t \\ &= N^*(\psi_{s,t})_*
A(\widetilde{1_B})_t , \endalign $$ so by the isomorphism theorem (3.6),
$A(\widetilde{1_B})_s= (\psi_{s,t})_* A(\widetilde{1_B})_t$. Hence $A
(\widetilde{1_B})_t = \tilde \beta_t$ for some constructible function
$\beta$ supported on $B$. The specialization theorem (3.7) yields $$\align
N^*A (1_{B}) &= (-1)^na_*N^*({1_{B}}) \\
 &=(-1)^na_*\lim_{t \downarrow 0} N^*(\widetilde{1_B})_t \ \ \ \ \text{by
(1)}\\ &=\lim_{t \downarrow 0} (-1)^na_*N^*(\widetilde{1_B})_t \\
&=\lim_{t \downarrow 0} N^*A(\widetilde{1_B})_t \\ &=\lim_{t \downarrow
0}N^*\tilde\beta_t\\ &= N^*(\beta), \endalign $$ so $\beta = A (1_B)$ by
(3.6). Thus $A(\widetilde{1_B})_t = (A(1_B))^{\widetilde{}}_t$. On the
other hand, $A (1_B) = D (1_B)$ by a direct computation, and
$(D(1_B))^{\widetilde{}}_t = D(\widetilde{1_B})_t$, which gives the claim
(3).

Therefore by (3.7) and (3.9), $$\align N^* A (1_\Delta) &= N^*A f_*(1_{B})
\\
 &=(-1)^na_*N^*f_* (1_{B}) \\
 &=(-1)^na_*\lim_{t \uparrow 1}N^*(\widetilde{1_B})_t\ \ \ \ \text{by (2)}
\\ &= \lim_{t \uparrow 1}(-1)^na_*N^*(\widetilde{1_B})_t \\ &= \lim_{t
\uparrow 1}N^*A(\widetilde{1_B})_t \\ &=\lim_{t \uparrow 1}
N^*D(\widetilde{1_B})_t\ \ \ \ \text{by (3)} \\ &= N^*f_*D (1_{B})\ \ \ \
\text{by (2)} \\ & = N^*D (1_\Delta), \endalign$$ from which we conclude
by (3.6) that $A (1_\Delta) = D (1_\Delta)$. \qed\enddemo

\proclaim{3.13 Corollary} Let $\alpha$ be a constructible function on the
oriented analytic manifold $M$. Then $$D_M(D_M\alpha) = \alpha. \tag{1}$$
If $\alpha$ has compact support, then $$\chi(D_M\alpha) = \chi(\alpha).$$
\endproclaim

\demo{Proof} By embedding $M$ in a Euclidean space, we may assume that the
support of $\alpha$ has codimension at least 2. Then (1) follows from
(3.12) and (3.6), and (2) follows from (3.12) and (2.10).\qed\enddemo

\remark{Remark} Corollary (3.13) has an elementary combinatorial proof,
which uses the triangulation theorem.\endremark

\proclaim{3.14 Proposition} Let $f:M\to N$ be a subanalytic map, and let
$\alpha$ be a constructible function on $M$ such that $f|\spt\alpha$ is
proper. Then $$D_N(f_*\alpha) = f_*(D_M \alpha).$$\endproclaim

\demo{Proof} As in the proof of (3.12), we may assume that $M$ and $N$ are
embedded in $\Bbb R^n$ as the ends of the mapping cylinder of $f$. We will
show that $D(f_*\alpha) = f_*(D\alpha),$ where $D=D_{\Bbb R^n}$. Let
$N^*=N^*_{\Bbb R^n}$. Using the notation of (3.9), we have $$\align
N^*D(f_*\alpha) &= (-1)^na_*N^*(f_*\alpha)\ \ \ \ \text{by (3.12)}\\ &=
(-1)^na_*\lim_{t \uparrow 1}N^*(\tilde\alpha_t)\ \ \ \ \text{by (3.7) and
(3.9)}\\ &= \lim_{t \uparrow 1}(-1)^na_*N^*(\tilde\alpha_t)\\ &= \lim_{t
\uparrow 1}N^*(D\tilde\alpha_t)\ \ \ \ \text{by (3.12)}\\ &=
N^*(f_*(D\alpha))\ \ \ \ \text{by (3.7) and (3.9),}\endalign$$ which
implies $D(f_*\alpha) = f_*(D\alpha)$ by (3.6).\qed\enddemo

\subhead 3.15 Mod 2 constructible functions\endsubhead Let $M$ be a real
analytic manifold. The function $\alpha: M \to \Bbb Z/2\Bbb Z$ is
(subanalytically) {\it constructible} if $\alpha^{-1} (0)$, and hence
$\alpha^{-1}(1)$, are subanalytic sets. The properties of $\Bbb Z/2 \Bbb
Z$-valued constructible functions are parallel to the properties of $\Bbb
Z$-valued constructible functions. The mod 2 Euler characteristic and the
mod 2 conormal cycle are defined as above: If $\alpha=\sum_in_i1_{X_i}$ is
a locally finite sum with $n_i\in\Bbb Z/2 \Bbb Z$ and $X_i$ closed
subanalytic, then $$\align \chi(\alpha)&=\sum_in_i\chi(X_i),\\ \Cal
N^*_M(\alpha)&=\sum_in_i\Cal N^*_M(X_i).\endalign$$ Here $\chi(X_i)$
denotes the mod 2 Euler characteristic of $X_i$, and $\Cal N^*_M(X_i)$ is
the mod 2 conormal cycle (2.13). Pushforward and duality of mod 2
constructible functions are defined as for $\Bbb Z$-valued constructible
functions, and all the results of this section are true for mod 2
constructible functions. In particular, the specialization theorem (3.7)
and the mapping cylinder proposition (3.9) are true mod 2; their proofs
use the {\it local} Gauss-Bonnet theorem. The $\Bbb Z/2 \Bbb Z$ versions
of (3.10)--(3.14) can be deduced from the $\Bbb Z$ versions by embedding
in Euclidean space and reduction mod 2.

\heading 4. Geometric Stiefel-Whitney classes \endheading

Sullivan \cite{S} discovered a local topological property of real analytic
varieties which implies that the combinatorial Stiefel chains (1.1) are
cycles, and hence that the Stiefel-Whitney homology classes are defined.
We show that Sullivan's local Euler characteristic condition for a
subanalytic set $X$ is equivalent to the antipodal symmetry of the
conormal cycle of $X$. This observation leads to a new definition of
Stiefel-Whitney classes. In order to prove the Deligne-Grothendieck
pushforward formula, we define the Stiefel-Whitney classes for certain
constructible functions on an analytic manifold.

The following definition is due to Sullivan \cite{S} \cite{A}.

\definition{4.1 Definition} A (mod 2) {\it Euler space} is a triangulable
topological space $X$ such that for all $x\in X$ the local Euler
characteristic $\chi_x(X) = \sum_i \rank H_i(X,X\setminus\{x\})$ is
odd.\enddefinition

If $\lk_x(K)$ is the link of $x$ in a triangulation $K$ of $X$ then a
neighborhood of $x$ in $X$ is homeomorphic to the cone on $\lk_x(K)$, so
$\chi_x(X)=1-\chi(\lk_x(K))$. Thus $X$ is an Euler space if and only if
$\chi(\lk_x(K))$ is even for every $x\in X$ and every triangulation $K$ of
$X$. A subanalytic set $X$ is an Euler space if and only if, for every
$x\in X$, the subanalytic link of $x$ in $X$ ({\it cf\.} \cite{CK}) has
even Euler characteristic.

We generalize Sullivan's definition as follows.

\definition{4.2 Definition} The mod 2 constructible function $\alpha$ on
the analytic manifold $M$ is an {\it Euler function} if $$D_M\alpha=
\alpha.$$\enddefinition

By the definition of the duality operator $D_M$, the subanalytic subset
$X$ of the analytic manifold $M$ is an Euler space if and only if the
characteristic function $1_X$ is an Euler function.

The following result is an immediate corollary of Theorem 3.12. Recall
that $a:S^*M\to S^*M$ is the antipodal map of the cotangent ray space of
$M$, and $\Cal N^*_M(\alpha)$ is the mod 2 conormal cycle of $\alpha$.

\proclaim{4.3 Theorem} The mod 2 constructible function $\alpha$ on $M$ is
an Euler function if and only if $$a_*\Cal N^*_M(\alpha)= \Cal
N^*_M(\alpha).\qed$$ \endproclaim

We will say that the mod 2 integral current $\Cal C$ in the cotangent ray
space $S^*M$ is {\it symmetric} if $a_*\Cal C= \Cal C $. Thus $\alpha$ is
Euler if and only if $\Cal N^*_M(\alpha)$ is symmetric, and the
subanalytic set $X$ is an Euler space if and only if the conormal cycle of
$X$ is symmetric.

Now we prove the following key result of Sullivan \cite{S}.

\proclaim{4.4 Theorem} Every real analytic set is an Euler space.
\endproclaim

\demo{Proof} Let $X$ be a real analytic set, which we may assume to be an
analytic variety in $\Bbb R^n$. Let $Z\subset\Bbb C^n$ be a
complexification of $X$. By \cite{F5} the conormal cycle of $Z$ in $\Bbb
C^n$ is invariant under multiplication by $\lambda$ for all complex
numbers $\lambda$ of modulus 1; in particular $a_*N^*_{\Bbb C^n}
(Z)=N^*_{\Bbb C^n}(Z)$. Therefore by (4.3) $Z$ is an Euler space. Now if
$x\in X$ and $\Sigma$ is a small sphere about $x$ in $\Bbb C^n$, then the
complex conjugation map $\tau:\Bbb C^n\to \Bbb C^n$ is an automorphism of
the link $Z\cap\Sigma$ with fixed-point set $X\cap\Sigma$. Therefore
$\chi(Z\cap\Sigma)\equiv\chi(X\cap\Sigma) \pmod{2}$. Thus
$\chi(X\cap\Sigma)$ is even.\qed\enddemo

The Stiefel-Whitney homology classes for a compactly supported subanalytic
Euler functions are constructed using the projectivized conormal cycle.
Let $\alpha$ be a mod 2 constructible function on the analytic
$n$-manifold $M$. Let $\pi:\Bbb PT^*M\to M$ be the projectivized cotangent
bundle of $M$. The conormal cycle $\Cal N_M^*(\alpha)$ is a closed
$(n-1)$-dimensional Legendrian integral current in the cotangent ray space
$S^*M$. Let $q:S^*M\to\Bbb PT^*M$ be the quotient map.

\definition{4.5 Definition} The {\it projectivized conormal cycle} of
$\alpha$ is the mod 2 current $$\Bbb PN^*_M(\alpha) = q_*(\Cal
N^*_M(\alpha)\ \LL\ U),$$ where $U$ is a subanalytic subset of $S^*(M)$
with $U\cap aU=\emptyset$ and $U\cup aU=S^*M.$\enddefinition

By (4.3), for $U$ as above, the conormal cycle satisfies $$a_*(\Cal
N_M^*(\alpha)\ \LL\ U) = a_*\Cal N^*_M(\alpha)\ \LL\ aU = \Cal
N^*_M(\alpha)\ \LL\ aU.$$ Therefore, $$\align a_*\partial(\Cal
N^*_M(\alpha)\ \LL\ U) &= \partial a_*(\Cal N^*_M(\alpha)\ \LL\ U) \\ &=
\partial(\Cal N^*(\alpha)\ \LL \ aU) \\ &= \partial (\Cal N^*_M(\alpha)\
\LL\ U) \\ \intertext {since $ \partial(\Cal N^*_M(\alpha)\ \LL\ U + \Cal
N^*_M(\alpha)\ \LL\ aU) =\partial \Cal N^*_M(\alpha)=0$, so} \partial \Bbb
P N^*_M(\alpha) &= q_* \partial (\Cal N^*_M(\alpha)\ \LL\ U) = 0.
\endalign $$

Furthermore, $\Cal N^*_M(\alpha)$ may be reconstructed from $\Bbb
PN^*_M(\alpha)$ as the fiber product $$\Cal N^*_M(\alpha)= \Bbb
PN^*_M(\alpha)\times_{\Cal B}\left(\lb +1\rb +\lb -1\rb\right),$$ where
$\Cal B$ is the bundle $q:S^*M \to \Bbb PT^*M$; this property is
sufficient to determine $\Bbb PN^*_M(\alpha)$, and shows that it is
well-defined, independent of the choice of the set $U$ above. Furthermore,
if $\spt\alpha$ is contained in the subanalytic set $X$, then $\spt \Bbb
PN^*_M(\alpha)\subset\pi^{-1}(X)$. Let $[\Bbb PN^*_M(\alpha)]$ be the mod
2 homology class of $\Bbb PN^*_M(\alpha)$ in $\pi^{-1}(X)$. Let
$\pi_X:\pi^{-1}(X)\to X$ be the restriction of $\pi:\Bbb P T^*M\to M$. We
abbreviate $\Bbb PN^*_M(\alpha)$ to $\Bbb PN^*(\alpha)$.

Let $\zeta_M \in H^1(\Bbb P T^*M;\Bbb Z/2\Bbb Z)$ be the mod 2 Euler class
(first Stiefel-Whitney class) of the tautological line bundle over $\Bbb
PT^*M$, and let $w^i(M) \in H^i(M;\Bbb Z/2\Bbb Z)$, $i=0,\dots,n$, be the
Stiefel-Whitney classes of the tangent bundle of $M$. Consider the
cohomology classes $$\gamma^k_M = \sum_i \zeta^{k-i}_M \smile w^i(M)$$ in
$H^k(\Bbb PT^*M;\Bbb Z/2\Bbb Z)$, $k\geq 0$, where $\smile$ denotes cup
product.

\definition{4.6 Definition} Let $\alpha$ be a subanalytic Euler function
on the analytic manifold $M$, with compact support contained in the
subanalytic set $X$. For each $i\geq 0$, the $i$th {\it Stiefel-Whitney
class} of $\alpha$ in $X$ is
 $$ (w_X)_i (\alpha) = (\pi_X)_*\left([\Bbb
PN^*(\alpha)]\frown\gamma^{n-i-1}_M\right)$$ in $H_i(X;\Bbb Z/2\Bbb Z),$
where $\frown$ denotes cap product. \enddefinition

We will abbreviate $(w_X)_i(\alpha)$ to $w_i(\alpha)$; the set $X$ will be
clear from context. This definition is parallel to Fu's definition of the
Chern-MacPherson homology classes of a complex analytic variety \cite{F5,
2.4}.

First we observe that the Stiefel-Whitney classes are additive.

\proclaim{4.7 Proposition} Let $\alpha$ and $\beta$ be constructible
functions on $M$ with compact support contained in the subanalytic set
$X$. Then, for all $i\geq 0$, $$w_i(\alpha + \beta) = w_i(\alpha) +
w_i(\beta).$$ \endproclaim

\demo{Proof} We have $$\align w_i(\alpha+\beta)&=(\pi_X)_*\left([\Bbb
PN^*(\alpha+\beta)]\frown \gamma^{n-i-1}_M\right) \\ &=
(\pi_X)_*\left([\Bbb PN^*(\alpha)+\Bbb PN^*(\beta)]\frown
\gamma^{n-i-1}_M\right) \\ &= (\pi_X)_*\left(\left[\Bbb
PN^*(\alpha)\right]\frown \gamma^{n-i-1}_M\right) +(\pi_X)_*\left([\Bbb
PN^*(\beta)]\frown \gamma^{n-i-1}_M\right) \\ &=
w_i(\alpha)+w_i(\beta).\qed\endalign$$\enddemo

Next we check that our definition agrees with the classical definition for
manifolds. If $E$ is a vector bundle of rank $r$ and $\zeta$ is the Euler
class of the tautological line bundle on $\Bbb PE$, then we have the
Wu-Hirsch relation \cite{Hu, 16.2.6} $$\sum_i\zeta^{r-i} \smile w^i(E) =
0. \tag{4.8}$$ If the base $B$ of $E$ is a compact manifold and $p:\Bbb
PE\to B$ is the projection, (4.8) implies that for all $i\geq 0$
$$p_*\zeta^{r+i-1}= \bar w^i(E),\tag{4.9}$$  where $p_*$ is the Gysin
homomorphism (fiber integration), $\bar w^0=1$ and $$\sum_i\bar
w^{k-i}(E)\smile w^i(E)=0$$ for all $k>0$.

In the following we let $w_\bull(\alpha)=w_0(\alpha)+w_1(\alpha)+\dots\in
H_*(X;\Bbb Z/2\Bbb Z)$ and $w^\bull(E)=w^0(E)+w^1(E)+\dots\in H^* (B;\Bbb
Z/2\Bbb Z)$ denote the total Stiefel-Whitney classes.

\proclaim{4.10 Proposition} If the subanalytic set $X$ in the analytic
manifold $M$ is a compact $C^1$ submanifold of dimension $d$, then
$w_i(1_X)$ is the Poincar\' e dual of the Stiefel-Whitney cohomology class
$w^{d-i}(X)$, $i=0,\dots,d$. \endproclaim

\demo{Proof} Since $X$ is a $C^1$ submanifold of $M$, the projectivized
conormal cycle of $X$ is the mod 2 integral current given by the
projectivized conormal bundle $\Bbb P(n^*X)$. Now $$\align w_\bull(1_X)
&=(\pi_X)_*\left([\Bbb PN^*(1_X)]\frown\left(\tsize\sum\zeta^i_M\smile
w^\bull(M)\right)\right)\\ &=(\pi_X)_*\left([\Bbb
PN^*(1_X)]\frown\tsize\sum\zeta^i_M\right)\frown w^\bull(M)\\
&=\left([X]\frown \bar w^\bull(n^*X)\right)\frown w^\bull(M) \ \ \ \
\text{by (4.9)}\\ &= [X]\frown\left(\bar w^\bull(n^*X)\smile
w^\bull(M)\right)\\ &= [X]\frown\left(\bar w^\bull(n^*X)\smile\left(
w^\bull(n^*X)\smile w^\bull(X)\right)\right)\\ &= [X]\frown
w^\bull(X),\endalign $$ as claimed. \qed \enddemo

Next we prove a special case of the fact that the Stiefel-Whitney classes
of the Euler function $\alpha$ on $X$ are independent of the embedding of
$X$.

\proclaim{4.11 Lemma} Let $M$, $M'$ be analytic manifolds, let
$\varphi:M\to M'$ be an analytic embedding, and let $\alpha$ be a
compactly supported Euler function on $M$. Then $\varphi_*w_\bull (\alpha)
= w_\bull (\varphi(\alpha))$.\endproclaim

\demo{Proof} Let $n=\dim M$, $n'=\dim M'$. Identifying $M$ with
$\varphi(M)$, the product formula (2.9) implies that the projectivized
conormal cycle $\Bbb PN^*_{M'}(\alpha)$ may be expressed as the join of
$\Bbb PN^*_{M}(\alpha)$ with the projectivized normal bundle of $M$. More
precisely, let $\widetilde P$ denote the blowup of $\Bbb P T^*M'|_M$ over
$\Bbb PT^* M \cup \PP (n^*M)$; {\it i.e.}, $$\align \widetilde P &= \{ \,
([\xi,\eta],[\xi'],[\eta'])\ |\ \xi \wedge \xi' = 0 = \eta \wedge \eta'
\,\} \tag{1} \\ &\subset (\PP T^*M'|_M) \times_M \PP T^*M \times_M \PP
(n^*M), \endalign $$ where the blowdown map $\sigma: \widetilde P \to \PP
T^*M'|_M$ is the projection onto the first factor. Then (2.9) implies that
the blowup cycle $\widetilde N = (\sigma^{-1})_*( \PP N^*_{M'} (X))$, a
closed mod 2 integral current on $\widetilde P$, may be expressed as the
iterated fiber product $$\widetilde N= \left( \PP N^*_M(\alpha)
\times_{\Cal B}\lb\PP^{n'-n-1}\rb\right) \times_{\Cal S} \lb\PP^1\rb,
\tag{2}$$ where $\Cal B$ is the pullback to $\PP T^*M$ of the
$\PP^{n'-n-1}$ bundle $\PP (n^*M) \to M$, and $\Cal S$ is the $\PP^1$
bundle over $\PP T^*M \times_M \PP (n^*M)$ given by the projection onto
the second and third factors in (1).

Next observe from (1) that the pullback via $\sigma$ of the tautological
line bundle over $\PP T^* M'|_M$ is the tautological bundle over $\Cal S$:
$$\sigma^* \Cal O_{\PP T^*M'\mid_M} (1) = \Cal O_{\Cal S} (1).\tag{3}$$ On
the other hand, as a bundle over $\PP T^* M \times_M \PP (n^* M)$, $$\Cal
S = \PP \bigl(\Cal O_{\PP T^*M}(1) \oplus \Cal O_{\PP (n^*M)}(1)\bigr),$$
so if we put $\zeta = w^1(\Cal O_{\PP T^*M'\mid_M}(1))$, $z = w^1 (\Cal
O_{\PP T^*M}(1))$, $\bar z = w^1(\Cal O_{\PP (n^*M)}(1))$, then with (3)
the fundamental relation (4.8) gives $$\align 0 &= (\sigma^*\zeta)^2 +
(\sigma^*\zeta) \cdot w^1( \Cal O_{\PP T^*M} (1) \oplus \Cal O_{\PP (n^*
M)}(1)) + w^2(\Cal O_{\PP T^*M} (1) \oplus \Cal O_{\PP (n^* M)}(1)) \\ &=
(\sigma^*\zeta)^2 + (\sigma^*\zeta)(z +\bar z) + z\bar z  \tag{4}
\endalign$$ by the Whitney sum formula. Multiplying (4) by successive
powers of $\sigma^*\zeta$ and substituting for $(\sigma^* \zeta)^2$ via
(4), we find inductively, for each positive integer $k$,
$$(\sigma^*\zeta)^{k+1} + (\sigma^*\zeta)\left(\sum_{i+j=k} z^i \bar
z^j\right) + z \bar z \left(\sum_{i +j = k-1} z^i \bar z^j \right)= 0.
\tag{5}$$ Let $w_\bull(\alpha)$ be the Stiefel-Whitney class of $\alpha$
as a function on $M$, and let ${w_{\bull}}'(\alpha)$ be the
Stiefel-Whitney class of $\alpha$ as a function on $M'$. We have the
following computation of ${w_{\bull}}'(\alpha)$. (For simplicity we use
$\cdot$ to denote both cup and cap products.) $$\align
{w_{\bull}}'(\alpha)&= \pi_*\left([\PP N^*_{M'}(\alpha)] \cdot
\left(\tsize\sum \zeta^k\right) \cdot w^\bull (T^*M'|_M ) \right)\\ &=
\pi_*\left([\widetilde N] \cdot \left(\tsize\sum (\sigma^* \zeta)^k\right)
\cdot w^\bull (T^*M'|_M )\right) \\ &= \pi_*\left([\widetilde N] \cdot
(\sigma^*\zeta + z \bar z)\cdot\left(\tsize\sum z^i
\right)\cdot\left(\tsize\sum \bar z^j\right) \cdot w^\bull (T^*M'|_M
)\right)\ \ \ \  \text{by (5)}\\ &= \pi_*\left(\left(\left[\PP N^*_M(X)
\times_{\Cal B} \lb\PP^{n'-n-1}\rb\right]\right) \cdot \left(\tsize\sum
z^i \right)\cdot\left(\tsize\sum \bar z^j\right) \cdot w^\bull (T^*M)
\cdot w^\bull(n^*M)\right),\endalign $$ by (2), since $\sigma^*\zeta$
pairs nontrivially with the fibers of $\Cal S$. Factoring the projection
$\pi$ through the projection $\mu$ of $\Cal B$, $${w_{\bull}}'(\alpha) =
\pi_*\left([\PP N^*_M(\alpha)] \cdot \left(\tsize\sum z^i\right)\cdot
w^\bull(T^*M) \cdot \mu_*\left(\tsize\sum \bar z^j\right) \cdot
w^\bull(n^*M)\right).$$ Now the identity (4.8) for the bundle $E = n^*M$
may be written $$\mu_*\left(\tsize\sum\bar z^j\right)\cdot w^\bull(n^*M) =
1,$$ whence $${w_{\bull}}'(\alpha) = \pi_*\left([\PP N^*_M(X) ] \cdot
\left(\tsize\sum z^i\right) \cdot w^\bull (T^*M)\right) = w_\bull
(\alpha),$$ as desired.\qed\enddemo

Now we prove the basic pushforward property for Stiefel-Whitney classes.

\proclaim{4.12 Theorem} Let $\alpha$ be a mod 2 constructible function on
the analytic manifold $M$, such that the support of $\alpha$ is compact
and contained in the subanalytic set $X$. Let $Y$ be a subanalytic subset
of the analytic manifold $N$, and let $f:X\to Y$ be a subanalytic map.
Then, for all $i\geq 0$, $$f_*w_i(\alpha) = w_i(f_*\alpha).$$ \endproclaim

\demo{Proof} By (4.11) we may assume $M=\Bbb R^m$ and $N=\Bbb R^m$. Let
$C$ be the mapping cylinder of $f$ (3.9). Let $\lambda:X\to C$ and
$\mu:Y\to C$ be the inclusion maps $$\align\lambda(x)&=(x,0,f(x)),\\
\mu(y)&=(0,1,y).\endalign$$ Then $\mu$ is a homotopy equivalence, with
homotopy inverse $\nu (x,t,y)=y$, and $\lambda$ is homotopic to $\mu\circ
f$, by the homotopy $h_t(x)= ((1-t)x,t,f(x))$. Therefore if we reduce the
relation (3.9)(3) modulo 2, the definition of Stiefel-Whitney classes
(4.6) gives $$\align \mu_*w_\bull(f_*\alpha)&=\lambda_*w_\bull(\alpha)\\
&=\mu_*f_*w_\bull(\alpha),\endalign$$ so
$w_\bull(f_*\alpha)=f_*w_\bull(\alpha)$.\qed \enddemo

\remark{Remark} Our proof of the pushforward formula for Stiefel-Whitney
homology classes is analogous to the original combinatorial proof sketched
by Sullivan using the mapping cylinder \cite{S, p\. 167}. \endremark

We prove a specialization formula for the Stiefel-Whitney classes.

\proclaim{4.13 Proposition} Let $\alpha$ be a mod 2 constructible function
on the analytic manifold $M$ such that $X=\spt\alpha$ is compact and has
positive codimension. Let $g:M \to \Bbb R$ be a subanalytic function. For
$t\in \Bbb R$, put $\alpha_t = \alpha\cdot 1_{X_t}$. Suppose that there
exists $\epsilon>0$ such that $\alpha_t$ is Euler for all $t\in
(0,\epsilon]$ and $\Cal N^*(\alpha)\ \LL\ (g\circ\pi)^{-1}(0,\epsilon]$ is
symmetric. Then the specialization $\phi_0\alpha$ is Euler, and
$$\lim_{t\downarrow 0} w_\bull (\alpha_t) = w_\bull(\phi_0\alpha)$$ in the
following sense. For every neighborhood $U$ of $X_0$, for $t>0$
sufficiently small, $X_t\subset U$ and $$j_{t*} w_\bull(\alpha_t) = j_{0*}
w_\bull (\phi_0\alpha),$$ where $j_s: X_s \to U$ is the inclusion map.
\endproclaim

\demo{Proof} Given $U$, choose  $h$ so that $X_{[0,h]} \subset U$. Since
$\Cal N^*(\alpha)\ \LL\ (g\circ\pi)^{-1}(0,\epsilon]$ is symmetric, the
mod 2 current $Z_h$ of (3.7) is antipodally symmetric, and hence it admits
a projectivization $\Bbb P Z_h$. Furthermore $\Cal N^*(\phi_0\alpha) =
\Cal N^*(\alpha_h)-\partial Z_h$ is antipodally symmetric, since
$\alpha_h$ is Euler. Let $\pi_U:\pi^{-1}(U)\to U$ be the restriction of
$\pi:\Bbb PT^*M \to M$. Now $$\align j_{h*}w_\bull (\alpha_h) - j_{0*}
w_\bull(\phi_0\alpha) &= (\pi_U)_*\left(\left([\Bbb PN^*(\alpha_h)] -
[\Bbb P N^*(\phi_0 \alpha)]\right) \cdot \gamma^\bull\right)\\ &=
(\pi_U)_*([\partial \Bbb P Z_h]\cdot \gamma^\bull)\\ &= 0, \endalign$$ as
desired. \qed \enddemo

\remark{Remark} A specialization formula for bivariant Stiefel-Whitney
classes of piecewise-linear spaces has been proved by Fulton and
MacPherson \cite{FM} ({\it cf\.} \cite{E}), but under the assumption that
the parameter map $g$ is piecewise-linear. The specialization formula for
Chern classes is due to Verdier \cite{V}; other proofs of Verdier's
formula have been given by Sabbah \cite{Sa} and Fu \cite{F2}.\endremark

The parallel between our definition of Stiefel-Whitney classes and Fu's
definition of Chern classes gives the following result.

\proclaim{4.14 Theorem} Let $a$ be a complex analytically constructible
$\Bbb Z$-valued function on the complex manifold $M$ such that $a$ has
compact support contained in the complex analytic variety $X$. Let
$\alpha$ be the mod 2 equivalence class of $a$. Then, for all $i\geq 0$,
$$\align w_{2i}(\alpha) &= r(c_i(a)),\\ w_{2i+1}(\alpha) &=0,\endalign$$
where $r:H_{2i}(X;\Bbb Z) \to H_{2i}(X;\Bbb Z/2\Bbb Z)$ denotes reduction
mod 2. \endproclaim

\demo{Proof} Let $n$ be the complex dimension of $M$. Denote the real and
complex projectivizations of $T^*M$ by $\RP T^*M$ and $\CP T^*M$
respectively. Thus $\RP T^*M$ is an $\RP^1$-bundle $\Cal H$ over $\CP
T^*M$ with projection $\mu$. The pullback of the tautological line bundle
$\Cal O_{\Bbb C} (1) $ over $\CP T^*M$  may be written $$ \mu^* \Cal
O_{\Bbb C}(1) = \Cal O_{\Bbb R}(1) \oplus \sqrt{-1}\Cal O_{\Bbb R}(1)
\cong \Cal O_{\Bbb R}(1) \oplus \Cal O_{\Bbb R}(1)$$ as a real vector
bundle, so if we put $\zeta = c^1(\Cal O_{\Bbb C}(1))$, $z = w^1(\Cal
O_{\Bbb R}(1))$, then $$\mu^*\zeta \equiv w^2(\Cal O_{\Bbb R}(1) \oplus
\Cal O_{\Bbb R}(1))\equiv z^2 \pmod{2}$$ by the Whitney sum formula. Let
$\RP N^*(\alpha)$ denote the projectivized conormal cycle of $\alpha$; it
is a $(2n-1)$-dimensional mod 2 integral current in $\RP T^*M$. Let $\CP
N^*(a)$ denote the complex projectivized conormal cycle of $a$ in the
sense of \cite{F5}; it is a $(2n-2)$-dimensional integral current in $\CP
T^*M$. These conormal cycles satisfy $$\RP N^*(\alpha) \equiv \CP N^*(a)
\times_{\Cal H} \lb\RP^1\rb \pmod{2}. $$ Therefore $$\align w_\bull
(\alpha) &= \pi_*\left([\RP N^*(\alpha)] \cdot \left(\tsize\sum
z^i\right)\cdot w^\bull(T^*M)\right) \\ &\equiv \pi_* \left(\left[\CP
N^*(a) \times_{\Cal H} \lb\RP^1\rb\right]
 \cdot (1+z) \cdot \mu^*\left(\tsize\sum \zeta^j\right)\cdot
c^\bull(T^*M)\right) \pmod{2}. \endalign $$ The fibers of $\Cal H$ pair
nontrivially with $z$, so if we factor the projection $\pi$ through $\mu$
then fiber integration gives $$\align w_\bull(\alpha) &\equiv
\pi_*\left([\CP N^*(a)] \cdot \left(\tsize\sum \zeta^i\right)\cdot
c^\bull(T^*M)\right) \pmod{2} \\
 &= c_\bull(a). \qed\endalign $$\enddemo

Finally we show that our definition (4.6) of Stiefel-Whitney classes
satisfies the Deligne-Grothendieck axioms.

If $M$ is an analytic manifold, let $E(M)$ be the group of Euler functions
on $M$. For $\alpha\in E(M)$, let $\spt\alpha$ denote the support of
$\alpha$. If $X$ is a subanalytic subset of the analytic manifold $M$, let
$$F(X)=\{\alpha\in E(M)\ |\ \spt\alpha\subset X,\ \spt\alpha\
\text{compact}\ \}.$$ If $Y$ is a subanalytic subset of the analytic
manifold $N$, and $f:X\to Y$ is a subanalytic map, then the pushforward
$f_*\alpha$ (3.8) is an Euler function by (3.14), and $f_*\alpha$ has
compact support contained in $Y$. Thus $f_*:F(X)\to F(Y)$, and if $g:Y\to
Z$, then $(g\circ f)_* = g_*\circ f_*$ by (3.10). Thus $F$ is a functor
from subanalytic sets to abelian groups, which we call the {\it Euler
constructible function functor}.

For each integer $i\geq 0$ let $H_i$ be the functor from subanalytic sets
to abelian groups given by $H_i(X)=H_i(X;\Bbb Z/2\Bbb Z)$, the $i$th mod 2
homology group of $X$. If $X$ is a compact $d$-dimensional analytic
manifold, let $[X]\in H_d(X;\Bbb Z/2\Bbb Z)$ be the fundamental class of
$X$, and for each $j\geq 0$ let $w^j(X)\in H^j(X;\Bbb Z/2\Bbb Z)$ be the
$j$th Stiefel-Whitney class of the tangent bundle of $X$. Now we state the
fundamental existence and uniqueness theorem for Stiefel-Whitney homology
classes. It is a generalization of Sullivan's theorem (1.1).

\proclaim{4.15 Theorem} For each $i\geq 0$ there exists a unique additive
natural transformation $w_i$ from $F$ to $H_i$ such that if $X$ is a
compact analytic $d$-manifold then $w_i(1_X)$ is Poincar\'e dual to
$w^{d-i}(X)$.\endproclaim

In other words, $w_i$ is a function which assigns to each Euler function
$\alpha$ with compact support contained in $X$ an $i$-dimensional mod 2
homology class $w_i(\alpha)$ of $X$ such that \roster \item
$w_i(\alpha+\beta)=w_i(\alpha)+w_i(\beta)$ for all $\alpha$, $\beta$,
\item $f_*w_i(\alpha)=w_i(f_*\alpha)$ for all $\alpha$ and all $f:X\to Y$,
\item $w_i(1_X)=[X]\frown w^{d-i}(X)$ for all compact analytic
$d$-manifolds $X$. \endroster These properties are analogous to the
Deligne-Grothendieck axioms for Chern classes of complex analytic
varieties \cite{M}. Fulton and MacPherson \cite{FM, Thm. 6A, p\. 64}
(corrected by Elhauoari \cite{E}) have proved a bivariant version of
(4.15), but only in the piecewise-linear category.

\demo{Proof}To prove existence we define $w_i$ as in (4.6). Additivity (1)
was proved in (4.7), pushforward (2) was proved in (4.12), and the
manifold property (3) was proved in (4.10).

Uniqueness follows from the representability of mod 2 homology by analytic
manifolds, as in the proof of Sullivan's theorem (1.1). Let $X\subset M$,
and let $\alpha$ be an Euler function on $M$ such that $\spt \alpha$ is
compact and contained in $X$. Let $k$ be the dimension of $\spt\alpha$.
Let $\Cal S$ be a subanalytic stratification of $\spt \alpha$. Since
$\alpha$ is Euler, each $(k-1)$-dimensional stratum of $\Cal S$ is
incident to an even number of $k$-dimensional strata. Thus the union of
the $k$-dimensional strata of $\Cal S$ is a mod 2 subanalytic cycle $Z$
(which represents $w_k(\alpha)$).

We claim that there exists a compact analytic $k$-manifold $V$ and a
subanalytic map $f:V\to\spt Z$ of odd degree. By the Thom representability
theorem \cite{T, Thm\. III.2}, applied to the fundamental class $[Z]\in
H_k(\spt Z;\Bbb Z/2\Bbb Z)$, there exists a compact $C^\infty$ manifold
$V$ and a continuous map $f:V\to \spt Z$ such that $f_*[V]=[Z]$. Since
every compact $C^\infty$ manifold has an analytic structure ({\it cf\.}
\cite{H, Thm\. 4.7.1}), we can assume $V$ is analytic. Finally we choose
subanalytic triangulations of $V$ and $\spt Z$, and then use simplicial
approximation to replace $f$ by a simplicial (hence subanalytic) map
homotopic to $f$.

Therefore if $\beta=\alpha-f_*(1_V)$ then $\beta$ is Euler and
$\dim(\spt\beta) <k$. It follows by induction on $k$ that there exists a
subanalytic map $\varphi:W\to X$ such that $W$ is a finite disjoint union
of compact analytic manifolds and $\varphi_*(1_W)=\alpha$. Thus
$w_i(\alpha)$ is uniquely determined by properties (1), (2) and (3).\qed
\enddemo

\heading 5. Polar cycles \endheading

Let $X$ be a compact subanalytic Euler space in $\Bbb R^n$, and let $A$ be
a linear subspace of $\Bbb R^n$. We say that $A$ is {\it general} for $X$
if the following condition holds. Let $A^*$ be the subspace of $\Bbb
R^{n*}$ corresponding to $A$ under the isomorphism $\Bbb R^n\cong\Bbb
R^{n*}$ given by the Euclidean metric. Now $\Bbb PN^*(X)\subset\Bbb
PT^*(\Bbb R^n)\cong \Bbb R^n\times\Bbb P(\Bbb R^{n*})$, and $A$ is general
for $X$ if there exists a subanalytic Whitney stratification $\Cal S$ of
$\spt\Bbb PN^*(X)$ such that $\Bbb R^n\times \Bbb P(A^*)$ is transverse to
the strata of $\Cal S$. For $k=1,\dots,n$ the set of $k$-planes $A$ in
$\Bbb R^n$ such that $A$ is general for $X$ is an open dense subset of the
Grassmannian of $k$-planes in $\Bbb R^n$.

\definition{5.1 Definition} The {\it polar cycle} of the compact
subanalytic Euler space $X$ with respect to a plane $A$ general for $X$,
with $1\leq \dim A \leq \dim X$, is the closed mod 2 current $$\sigma(X,A)
= \pi_*\bigl(\Bbb PN^*(X) \cdot (\Bbb R^n\times \Bbb P(A^*)) \bigr),$$
where $\pi:\Bbb PT^*\Bbb R^n\to \Bbb R^n$ is the projection.
\enddefinition

Here $\cdot$ denotes intersection of mod 2 currents. If $\dim A=i+1$, the
homology class of the polar cycle $\sigma(X,A)$ is the Stiefel-Whitney
class $w_i(X)$. The support of $\sigma(X,A)$ is contained in the set
$\{x\in X\ |\ \exists\, \xi\in A^*,\ (x,\xi)\in \text{spt}N^*(X) \}$, the
``singularity locus'' of the orthogonal projection $X\to A$.

In fact $\sigma(X,A)$ is the generalization to subanalytic sets of the
simplicial singularity cycle of Banchoff and McCrory \cite{B} \cite{Mc}.
Let $X\subset \Bbb R^n$ be a simplicial Euler space; {\it i.e.}, $X$ is an
Euler space together with the structure of a linear simplicial complex $K$
in $\Bbb R^n$: $X=|K|$. Let $f:X\to\Bbb R^{i+1}$ be a map which is linear
on each simplex of $K$. We say that $f$ is {\it nondegenerate} if, for
every vertex $v$ of $K$, and for every set of distinct vertices
$v_0,\dots,v_{i+1}$ of the simplicial star of $v$, the affine span of the
points $f(v_0),\dots,f(v_{i+1})$ is $\Bbb R^{i+1}$. The {\it Euler
singularity} of $f$ is the mod 2 simplicial $i$-chain $\Sigma(f)$ of $K$
defined as follows ({\it cf\.} \cite{Mc, p\. 373}). Let $S$ be an
$i$-simplex of $K$, let $P$ be the affine hyperplane of $\Bbb R^{i+1}$
spanned by $f(S)$, and let $P^+$, $P^-$ be the closures of the two
components of $\Bbb R^{i+1}\backslash P$. Let $L$ be the simplicial link
of $S$ in $K$, and let $L^\pm=L\cap f^{-1} (P^\pm)$. Then $\chi(L^+)\equiv
\chi(L^-) \pmod{2}$, and the coefficient of $S$ in the chain $\Sigma(f)$
is $1-\chi(L^+)\pmod{2}$.

\proclaim{5.2 Proposition} Let $X$ be a simplicial Euler space in $\Bbb
R^n$, let $A$ be a linear subspace of $\Bbb R^n$, and let $p:\Bbb R^n\to
A$ be the orthogonal projection. If $p|X$ is nondegenerate then $A$ is
general for $X$, and $$\sigma(X,A) = \Sigma(p|X).$$ \endproclaim

\demo{Proof} Let $N\subset T(\Bbb R^n)$ be the set of normal vectors of
$X$; {\it i.e.}, $(x,u)\in N$ if and only if $(x,\xi_u)\in \spt N^*(X)$,
where $\xi_u(v)=\langle u/|u|,v \rangle$. The set $N$ is a conical subset
of $T(\Bbb R^n)$; in other words, if $(x,u)\in N$ and $t>0$, then
$(x,tu)\in N$. The plane $A$ is general for $X$ if and only if there is a
conical subanalytic Whitney stratification $\Cal G$ of $N$ such that $\Bbb
R^n\times A$ is transverse to the strata of $\Cal G$.

Let $X=|K|$. We will show that there is a canonical conical stratification
$\Cal G$ of $N$ associated to the triangulation $K$, and that $\Bbb
R^n\times A$ is transverse to the strata of $\Cal G$ if $p$ is
nondegenerate.

If $S$ is a (closed) $k$-simplex of $K$, let $S^o$ be the interior of $S$
in the affine $k$-plane spanned by $S$, and let $S^\bot$ be the
$(n-k)$-dimensional linear subspace of $\Bbb R^n$ orthogonal to $S$. Let
$\Cal F(S)$ be the conical stratification of $S^o\times S^\bot$ defined as
follows. Let $\Cal T$ be the set of all $(k+1)$-simplices of $K$ which
have $S$ as a face. For each $T\in\Cal T$, $T^\bot$ is a hyperplane of
$S^\bot$. We let $\Cal E(S)$ be the stratification of $S^\bot$ induced by
this collection of hyperplanes. More precisely, if $v$, $w\in S^\bot$, let
$[v,w]$ be the line segment between $v$ and $w$. Then $v$ and $w$ are in
the same stratum of $\Cal E(S)$ if for all $T\in\Cal T$ either
$[v,w]\subset T^\bot$ or $[v,w]\cap T^\bot = \emptyset$. Let $\Cal
F(S)=\{S^o\times E\ |\ E\in\Cal E(S)\}$.

Let $N'=\cup_S (S^o\times S^\bot)$, and consider the stratification $\Cal
F$ of $N'$ given by $\Cal F=\cup_S\Cal F(S)$. Now $N'\subset N$, and it
follows from the multiplicity formula (2.3) that $N$ is a union of strata
of $\Cal F$. We define the stratification $\Cal G$ of $N$ by $\Cal G=
\{F\in\Cal F\ |\ F\subset N\}$.

Now $p|X:X\to A$ is nondegenerate if and only if for all $S\in K$, the
rank of $p$ restricted to the affine span of $S$ is $\min\{\dim S, \dim
A\}$. This is equivalent to the statement that $A$ is transverse to
$S^\bot$. Therefore if $p|X$ is nondegenerate then $\Bbb R^n\times A$ is
transverse to the strata of $\Cal F$, and so $\Bbb R^n\times A$ is {\it a
fortiori} transverse to the strata of $\Cal G$.

It follows from the definition of the polar cycle $\sigma(X,A)$ and the
multiplicity formula (2.3) that the multiplicity of the cycle
$\sigma(X,A)$ is constant on the interior of each $i$-simplex $S$ of $K$.
We proceed to give a formula for this multiplicity $\iota_S$.

If $T$, $U$ are simplices of $K$, we write $T<U$ if $T$ is a proper face
of $U$. Given $S\in K$, the following subcomplexes of $K$ are associated
to $S$: the boundary of $S$, $\partial S =\{T\ |\ T<S\}$; the star of $S$,
$\text{St}(S) = \{T\ |\ \exists U,\ T\leq U\geq S\}$; the boundary of the
star of $S$, $\partial \text{St}(S)= \{T\in \text{St}(S)\ |\ T\ngeq S\}$;
the link of $S$, $\text{Lk}(S) = \{T\in\partial\text{St}(S)\ |\ \nexists
U,\ T\geq U\leq S\}$. For any subcomplex $L$ of $K$, we let $|L|$ denote
the union of the simplices of $L$.

Let $S$ be an $i$-simplex of $K$. For notational simplicity assume that
the origin $0$ is in the interior of $S$. Choose $\epsilon>0$ so that
$\epsilon<d(0,|\partial\text{St}(S)|)$. Let $B_\epsilon$ be the sphere of
radius $\epsilon$ in $\Bbb R^n$ centered at $0$. Let $u$ be a unit vector
in the line $S^\bot\cap A$, and let $\xi(v)=\langle u,v \rangle$. Choose
$\delta>0$ so that $\delta<\epsilon$ and $\delta<|\xi(y)|$ for all $y$
such that $y\in T\cap S^\bot\cap \partial B_\epsilon$ for some
$(i+1)$-simplex $T$ with $S<T$. Then, by (2.2) and the product formula
(2.9), $$\align \iota_S &= \chi\bigl(X\cap S^\bot \cap B_\epsilon\cap
\xi^{-1}[-h,\infty)\bigr) \bigr|^{h=+\delta}_{h=-\delta}\\
&=\chi(C)-\chi(\ell^+),\endalign$$ where $$C=X\cap S^\bot \cap
B_\epsilon\cap \xi^{-1}[-\delta,\delta],$$ $$\ell^+=X\cap S^\bot \cap
B_\epsilon\cap \xi^{-1}(\delta).$$

Now $S^\bot\cap |\text{St}(S)|$ is the cone from $0$ over the set
$S^\bot\cap |\partial\text{St}(S)|$, which is homeomorphic to
$L=|\text{Lk}(S)|$. Thus $C$ is a cone, since it is the intersection of
the convex neighborhood $S^\bot \cap B_\epsilon\cap
\xi^{-1}[-\delta,\delta)$ of $0$ and the cone $S^\bot\cap |\text{St}(S)|$.
So $\chi(C)=1$, and $$\iota_S = 1-\chi(\ell^+).$$ By radial projection
$\ell^+$ is homeomorphic to $X\cap S^\bot \cap \partial B_\epsilon\cap
\xi^{-1}[\delta,\infty)$. By the choice of $\delta$, $$X\cap S^\bot \cap
\partial B_\epsilon\cap \xi^{-1}[0,\delta]\cong \bigl(X\cap S^\bot \cap
\partial B_\epsilon\cap \xi^{-1}(\delta)\bigl)\times[0,\delta].$$
Therefore $\ell^+\cong X\cap S^\bot \cap \partial B_\epsilon\cap
\xi^{-1}[0,\infty)$, which is homeomorphic by radial projection to
$S^\bot\cap |\partial\text{St}(S)|\cap \xi^{-1}[0,\infty)\cong
|\text{Lk}(S)|\cap \xi^{-1}[0,\infty)= L^+.$ Therefore
$$\iota_S=1-\chi(L^+),$$ which is the multiplicity of $\Sigma(p|X)$ along
$S$. Thus $\sigma(X,A)=\Sigma(p|X)$.\qed \enddemo

\remark{Remark} The Banchoff-McCrory description of the Euler singularity
cycle of a nondegenerate simplexwise-linear map can be generalized to
stratified maps, and the preceeding proposition can be generalized to this
context using stratified Morse theory \cite{GM}.\endremark

Note that if $X\subset \Bbb R^n$ and $f:X\to\Bbb R^{i+1}$ is a
nondegenerate simplexwise-linear map, then $\Sigma(f)=\Sigma(p|\Gamma)$,
where $p:\Bbb R^n\times\Bbb R^{i+1}\to \Bbb R^{i+1}$ is the projection and
$\Gamma$ is the graph of $f$.

Following Banchoff and McCrory, if $X=|K|$ we define a map $f_i:X\to \Bbb
R^{i+1}$, nondegenerate and linear on the simplices of the barycentric
subdivision $K^\prime$, by setting $$f_i(b\Delta^k) =
(k,k^2,\dots,k^{i+1}),$$ where $b\Delta^k$ is the barycenter of the
$k$-simplex $\Delta^k$. Recall that the $i$th Stiefel chain $s_i(K)$ is
the sum of all the $i$-simplices in the barycentric subdivision of $K$.

\proclaim{5.3 Proposition} $\Sigma(f_i) = s_i(K).$ \endproclaim

\demo{Proof} The proof given for manifolds in \cite{B, p.\ 345} goes
through for Euler spaces with only one change: The link of a simplex is
not necessarily a sphere, but it has even Euler characteristic.\qed
\enddemo

This result was used by Banchoff and McCrory to prove that if $X$ is a
simplicial Euler space, then the Euler singularity cycle $\Sigma(f)$ of a
nondegenerate simplexwise-linear map $f:X\to\Bbb R^{i+1}$ represents the
Stiefel-Whitney homology class $w_i(X)$. Prior to their work the only
definition of Stiefel-Whitney homology classes for singular spaces was the
combinatorial formula. Our geometric construction of Stiefel-Whitney
homology classes and verification of the axioms, together with (5.2) and
(5.3), gives a proof of the combinatorial formula for these classes. This
proof is new even for manifolds.

\remark{Remark} The present approach may also be applied to the smooth
case, {\it i.e\.} to embedded curvilinear ($C^1$) simplicial complexes.
Just as in section 4, such a complex $X$ in a manifold $M$ is an Euler
space (dropping the subanalyticity condition) if and only if its conormal
cycle $\Cal N^*_M(X)$ is antipodally symmetric. Therefore the
Stiefel-Whitney classes $w_i (X)$ may again be constructed by the formula
(4.6). That these classes do not depend on the (piecewise smooth)
embedding of $X$ follows from the arguments of (3.7) (specialization),
(4.11) (independence of analytic embedding), and (4.12) (pushforward). The
key points which require clarification here are the mass bounds on the
current $Z$ in the specialization formula for the mapping cylinder of a
piecewise smooth homeomorphism, and the existence of the blowup
$\widetilde N$ in the proof of (4.11). Both of these results follow from
the hypothesis that the maps in question are diffeomorphisms piecewise,
and therefore all slicing procedures are uniformly transverse to some
compatible triangulation. In particular, the expression for the
Stiefel-Whitney classes in terms of polar cycles yields the classical
combinatorial formula for the Stiefel-Whitney classes of a $C^1$
triangulated manifold \cite{HT}.\endremark

\enddocument